\documentclass[11pt,a4paper]{article}
\usepackage[left=2.8cm, right=2.8cm, top=3cm, bottom=3cm]{geometry}
\usepackage{graphicx}
\usepackage[T1]{fontenc}
\usepackage[utf8]{inputenc}
\usepackage[overload]{empheq}
\usepackage{setspace}
\usepackage{amsmath}
\usepackage{amssymb}
\usepackage{amsfonts}
\usepackage{esint}
\usepackage{relsize}
\usepackage{bm}
\usepackage[dvipsnames]{xcolor}
\usepackage{enumitem}
\usepackage{url}
\usepackage{pdfpages}
\usepackage{wrapfig}

\makeatletter
\def\patchsect#1\let\@svsec\@empty{#1\def\@svsec{\leavevmode\kern1sp\relax}}
\let\old@sect\@sect
\def\@sect{\expandafter\patchsect\old@sect}
\makeatother

\usepackage{amsthm}
\theoremstyle{definition}

\theoremstyle{plain}
\newtheorem{defi}{Definition}[section]
\newtheorem{lemma}[defi]{Lemma}
\newtheorem{theorem}[defi]{Theorem}
\newtheorem{prop}[defi]{Proposition}

\newtheorem{remark}[defi]{Remark}

\usepackage{titlesec}
\titleformat{\section}
  {\centering\normalfont\large\bfseries}{\thesection.}{1em}{}
\titleformat{\subsection}
  {\normalfont\large\bfseries}{\thesubsection.}{1em}{}
\titleformat{\subsubsection}
  {\normalfont\normalsize\itshape\color{MidnightBlue}}{\thesubsubsection.}{1em}{}

\usepackage{hyperref}
 \hypersetup{
    colorlinks,%
    citecolor=red,%
    filecolor=black,%
    linkcolor=black,%
    urlcolor=black
}

\usepackage{float}

\usepackage{subcaption}
\usepackage{icomma}
\usepackage{appendix}
\usepackage[backend=bibtex, backref=true, style=numeric,maxbibnames=99, url=false, isbn=false]{biblatex}
\DeclareMathOperator*{\Id}{Id}
\DeclareMathOperator*{\Tr}{Tr}
\DeclareMathOperator*{\Diag}{Diag}
\newcommand{\jap}[1]{\langle #1 \rangle}

\title{On the monotonicity of the entropy production in the Landau-Maxwell equation}
\date{\today}
\author{Côme Tabary\footnote{\texttt{tabary@imj-prg.fr}. Université Paris Cité and Sorbonne Université, CNRS, IMJ-PRG, F-75013 Paris, France and DMA, ENS, Université PSL, CNRS, 75005 Paris, France. \texttt{ORCID:0009-0005-6651-8640}. }}
\begin{document}

\maketitle

\begin{abstract}
\noindent We study the homogeneous Landau equation with Maxwell molecules and prove that the entropy production is non-increasing when the directional temperatures are well-distributed. It implies that for any initial condition with finite mass and energy, the entropy production is guaranteed to be non-increasing after a certain explicit time. We also obtain exponential decay of the entropy production. We provide an explicit formula for the derivative of the entropy production which may be a starting point for future progress towards full monotonicity.
\end{abstract}



\section{Introduction}

\subsection{Background and main result}

We consider in this work the homogeneous Landau equation in $d\geq 2$ space dimensions
\begin{equation}
\label{eq:landau}
    \partial_t f_t(v) = \nabla_v \cdot \int_{\mathbb{R}^d} \alpha(\vert v-w \vert) a (v - w) (\nabla_v - \nabla_w) f_t(v)f_t(w) dw,
\end{equation}
where the unknown $f_t:\mathbb{R}^d\rightarrow\mathbb{R}_+$ is the time-dependent distribution of velocities in a plasma.
The matrix-valued function
$$a(z):=\frac{1}{d-1}\left[\vert z \vert^2 \Id - z \otimes z\right]$$
is the projection on $z^{\perp}$ up to a factor $(d-1)^{-1}\vert z \vert^2$. The $d-1$ in the definition of $a$ is unusual but simplifies formulas, and amounts to a simple change of time. The non-negative function $\alpha$ is called the interaction potential, and is dictated by the nature of the interactions between the plasma particles.

The Landau equation first originated as a meaningful alternative to the Boltzmann equation in the case of 3D Coulomb interactions, where $\alpha(r)=r^{-3}$ —for which the latter equation fails to make sense. The Landau and Boltzmann equations are thus intimately linked, as the former can be obtained as a limit of the latter by the so-called \textit{grazing collisions limit} (see \cite{LifsicPitaevskij2008, Villani2002} for more details on kinetic models).

No matter the choice of $\alpha$, the Landau equation preserves mass, momentum and energy:
\begin{align*}
    \frac{d}{dt}\int_{\mathbb{R}^d} f_t(v)dv=\frac{d}{dt}\int_{\mathbb{R}^d} vf_t(v)dv=\frac{d}{dt}\int_{\mathbb{R}^d} \vert v\vert^2f_t(v)dv=0.
\end{align*}
Moreover, Boltzmann's H-theorem holds for the Landau equation, stating that the opposite of the entropy
$$H(f):=\int_{\mathbb{R}^d} f(v) \log f(v) dv$$
is non-increasing along solutions of \eqref{eq:landau}. A straightforward computation using the symmetry of the equation indeed shows that$\frac{d}{dt}H(f_t)=-D(f_t)$
where $D$ is given by
\begin{equation}
    \label{eq:D}
    D(f):=\frac{1}{2}\iint_{\mathbb{R}^{2d}} \alpha(\vert v-w\vert)a(v-w):\left[(\nabla_v-\nabla_w)\log\left(f(v)f(w)\right)\right]^{\otimes 2}f(v)f(w)dvdw \geq 0.
\end{equation}
We will refer to $D$ as the \textit{entropy production} functional. A natural follow-up question is: can we say anything on the variations of $D$? This question was first asked by McKean in his famous 1966 work \cite{McKean1966} on a one-dimensional toy model for the Boltzmann equation, known as Kac's model. Therein, he conjectured that 
for the solution $f=(f_t)_{t\geq0}$ of Kac's equation,
\begin{align}
\label{eq:mckean2}
    \frac{d}{dt}D(f_t)\leq 0,
\end{align}
and even suggested the \textit{complete monotonicity} of $H$, that is
\begin{align}
\label{eq:compmonH}
    (-1)^{n+1}\frac{d^n}{dt^n}H(f_t)\leq 0
\end{align}
for all $n\geq 1$. Shortly after, Harris \cite{Harris1967} proved the complete monotonicity for a further simplified toy model with discrete velocities. The conjecture then remained untouched, until 1982, when Olaussen \cite{Olaussen1982} and Lieb \cite{Lieb1982} independently disproved it, by studying the Bobylev-Krook-Wu particular solutions of the Boltzmann equation \cite{Bobylev1975a,KrookWu1976}. However, Lieb's abstract argument does not indicate for which $n$ property \eqref{eq:compmonH} fails, and Olaussen's proof shows that it fails for $n\geq 102$. Furthermore, numerical simulations \cite{ZiffMerajverStell1981} indicate that the complete monotonicity of $H$ holds for $n\leq 30$ (still for the particular Bobylev-Krook-Wu solutions).

In view of theses results, one might be tempted to think that McKean's conjecture may still be true for small values of $n$, especially for $n=2$. However, recently, Silvestre \cite{Silvestre2026} provided a counter-example to the monotonicity of entropy production in the case of the Boltzmann equation. Albeit with a non-physical collision kernel, the mechanism in the proof seems to indicate that monotonicity could also fail in physical cases. Yet, the counter-example does not seem to transfer to the Landau setting, suggesting that the conjecture \eqref{eq:mckean2} might still hold for certain Landau models.

In this work, we investigate the conjecture \eqref{eq:mckean2} for the Landau equation with Maxwell molecules, which corresponds to $\alpha=1$. While Coulomb interactions are the most physically relevant choice for $\alpha$, the simplified case of Maxwell molecules provides an easier-to-investigate model while still preserving some core properties of the equation. For the Boltzmann and Landau equations, the Maxwellian setting allows for much more explicit computations which can provide insight into the more physical settings: for instance, the monotonicity of the Fisher information along solutions of \eqref{eq:landau} was known to hold for Maxwell molecules \cite{Toscani1998,Villani2000, Villani1998a} more than two decades before the recent impressive extension \cite{GuillenSilvestre2023, ImbertSilvestreVillani2024} to other potentials. It thus seems reasonable to first study McKean's conjecture for Maxwell molecules before turning to any other potential $\alpha$.

Our main result shows monotonicity of the entropy production after a certain amount of time, relying on an explicit formula for the derivative of $D$. We show that this monotonicity in fact occurs as soon as the directional temperatures of the solution are sufficiently well-distributed in all directions. The knowledge of the exact evolution of the directional temperatures in the Maxwellian setting allows us to provide an explicit time after which monotonicity is guaranteed. We refer to \cite{Villani1998b} for the well-posedness of the Landau equation with Maxwell molecules for initial data with finite mass and energy. Up to rescaling the variable and choosing a suitable basis, we can always suppose that a function $f_0$ with finite mass and energy satisfies the normalization
 \begin{align}
    \label{eq:normalization}
        \int_{\mathbb{R}^d} f_0 dv=1, &&\! \int_{\mathbb{R}^d} v f_0 dv=0,&& \! \int_{\mathbb{R}^d} (v\otimes v) f_0 dv=T_0=\Diag(T_{0,i})_{1\leq i \leq d},&&\!  \sum_{i=1}^d T_{0,i}=d.
    \end{align}

\begin{theorem}
\label{thm:main}
    Let $f_0\geq 0$ with finite mass and energy satisfy the normalization \eqref{eq:normalization}
    and let $T_{0,\max}=\max_{1\leq i \leq d} (T_{0,i})\geq 1$. Consider $f=(f_t)_{t\geq 0}$ the solution to the Landau equation \eqref{eq:landau} with Maxwell molecules $\alpha=1$ and initial condition $f_0$.
    Then, there exists $t_0 \geq 0$ depending only on $d$ and $T_{0,\max}$ such that for all $t\geq t_0$, $D(f_t)$ is monotone non-increasing.
    We can further take 
    $$t_0=\frac{d-1}{4d}\log\left(\frac{T_{0,\max}-1}{C_d}\right)$$
    or $t_0=0$ if the above quantity is negative. The constant $C_d$ is given by
    \begin{align*}
       C_d= \begin{cases}
            \frac{13-\sqrt{133}}{18}\approx 0.08 & \text{if }d=2 \\
            \frac{20-2\sqrt{70}}{15}\approx 0.2 & \text{if }d=3 \\
         (d-1)^2\left((6d-3)+\frac{d(d+1)}{\eta(d-1)}\right)^{-1}\approx_{d\rightarrow \infty} 0.15d & \text{for }d\geq 4,
        \end{cases}
    \end{align*}
    where $\eta=\frac{9+\sqrt{129}}{12}\approx 1.69$.
\end{theorem}
\begin{remark}
    In dimension $3$, $T_{0,\max}\leq 3$ so the bound $t_0\leq \log(2/C_d)/6\leq 0.4$ always hold.
\end{remark}
Let us comment on this result. For radial data, the Landau-Maxwell equation reduces to the Fokker-Planck equation, and monotonicity of entropy production (which is then simply the Fisher information) is guaranteed. Theorem~\ref{thm:main} extends this monotonicity for general function that are \textit{not too anistropic}, where the anisotropy is quantified by how much the largest directional temperatures differs from $1$. We can interpret $[0,t_0]$ as some \textit{thermalization period} after which $f_t$ has become close enough to a radial function for monotonicity to hold. We point out that the constant $C_d$ and hence the time $t_0$ can certainly be improved on by technical optimizations of the proof, but we do not expect to cover all initial data by such means. We do believe that monotonicity actually holds for any initial data, but that it requires to better understand the favorable behaviour of some terms in the explicit formula for $\frac{d}{dt}D$, that the current proof ignores. This is further discussed in Remark~\ref{rem:fullmono}. 

Along the Fokker-Planck equation, corresponding to the radial solutions of the Landau-Maxwell equation, not only does the entropy production / Fisher information decrease, it actually decays exponentially with rate $2$. For the general solutions of the Landau-Maxwell equation, we can retrieve exponential decay of the entropy production with a rate arbitrary close to $2$:
\begin{theorem}
\label{thm:side}
     Let $f_0\geq 0$ with finite mass and energy satisfying the normalization \eqref{eq:normalization}, and let $f=(f_t)_{t\geq 0}$ be the solution to the Landau equation \eqref{eq:landau} with Maxwell molecules. Then, for any $0< \alpha< 2$, there exists $t_1>0$ depending only on $d$ and $\alpha$ such that for all $t\geq t_1$:
    $$D(f_t)\leq D(f_{t_1})e^{-\alpha (t-t_1)}.$$   
\end{theorem}
The proof of this theorem will be a simple combination of our explicit formula for $\frac{d}{dt}D$ and of the Bakry-Emery $\Gamma^2$-criterion on the sphere \cite{BakryGentilLedoux2014}. The rate $2^-$ is the same as obtained for the already known exponential decay of the relative entropy \cite{Villani1998b,DesvillettesVillani2000}.
    
\subsection{Literature review and discussion}

We discuss below several topics that motivated the present work, or that are connected to it.
\\

\noindent
\textbf{McKean's conjectures.} In his seminal 1966 paper \cite{McKean1966}, McKean showed the convergence of a one-dimensional $N$ particle system initially proposed by Kac \cite{Kac1956} towards the associated Boltzmann-like equation as $N\rightarrow+\infty$, a property now called \textit{propagation of chaos}. At the center of his proof was the monotonicity of the Fisher information (called  in his work \textit{Linnik's functional}) $
    I(f)=\int f(v) \vert \nabla \log f(v)\vert^2dv$
 for Kac's model. The Fisher information being the derivative of the entropy along the heat flow (in our notation, $\frac{d}{dt}H(f_t)=-i(f_t)$ for solutions of $\partial_t f_t =\Delta f_t$), McKean is lead to question the behavior of higher order time-derivatives of $H$. Let
\begin{equation*}
H^{(n)}(f):=(-1)^n\left.\frac{d^n}{dt^n}H(f_t)\right\vert_{t=0}
\end{equation*}
be the functional obtained by differentiating $n$ times the entropy along the heat flow with a probability density $f$ as initial condition. McKean conjectures that $H^{(n)}$ is minimized by the Maxwellian distribution with the same momentum and energy as $f$, which would imply in particular that the $H^{(n)}$ are all non-increasing along the heat flow: in other terms, $t\mapsto H(f_t)$ would be completely monotone. After several partial positive results depending on the dimension and convexity properties of $f$ \cite{FanChengYanlinGeng2015,Toscani2015,GuoYuanGao2022,ZhangAnantharamGeng2018,Wang2024}, an explicit counter-example in dimension $1$ was recently found by Gu and Sellke \cite{GuSellke2026}, where the conjecture fails at the fifth derivative. For a review of conjectures on the connections between entropy and the heat flow, we refer to \cite{Ledoux2022}.

McKean also conjectures that all the $H^{(n)}$ are monotone non-increasing along the solutions of Kac's Boltzmann-like equation, and that in some sense they are the only functionals of the form
$\int h^{(n)}(f,\partial_vf,...,\partial_v^nf)dv$
to do so (remark that $h^{(n)}$ has no explicit dependence on $v$, so this conjecture does not rule out monotonicity of $D$!). To our knowledge, no progress has been made on this conjecture past $n=0$ (entropy) and $n=1$ (Fisher information), but in view of the results on the heat flow, it would be very surprising for the conjecture to hold for all $n$.

Finally, McKean introduces the conjecture studied in this work, on the successive derivatives of the entropy along the Boltzmann or Landau equation, which lead to the previously discussed works \cite{Harris1967,ZiffMerajverStell1981,Olaussen1982,Lieb1982,Silvestre2026}. We simply add that none of these works treat the original Kac's $1$-dimensional Boltzmann model that is the focus of McKean's work: it is not known to us whether the counter example of Silvestre transfers to Kac's model.
\\

\noindent
\textbf{Maxwell molecules.} The Maxwellian setting $\alpha=1$ for the Landau equation stands as an intermediate between hard potentials $\alpha(r)=r^\gamma$ with $\gamma\in(0,1]$ and soft potentials $\gamma\in[-d,0)$. The interest of the Maxwellian case lies in the existence of an explicit expression for the matrix
$$\int_{\mathbb{R}^{d}} \alpha(\vert v-w\vert)a(v-w) f_t(w)dw ,$$
once the initial condition has been properly normalized. This explicit formula deletes the non-linearity of the equation (see Lemma~\ref{lem:linearlandau}). There is also a reliable expression for the Fourier transform of the equation by Bobylev \cite{Bobylev1975}, later used by Pulvirenti and Toscani \cite{PulvirentiToscani1996} to provide a direct proof of Tanaka's contractivity in Wasserstein distance \cite{Tanaka1978}. The Fisher information was known to decrease, since it did for the Boltzmann equation with Maxwell molecules and the result could be passed to the grazing collisions limit \cite{Toscani1998,Villani1998a}. We refer to the two works by Villani \cite{Villani1998b} for a complete wellposedness theory supplemented by many qualitative properties and convergence to equilibrium results, and \cite{Villani2000} for the first direct proof of monotonicity of the Fisher information.

As stated before, the Landau-Maxwell equation has often been considered as a toy model for the general setting, because the absence of $\alpha$ makes explicit computations possible. Recently, Caja-Lopez, Delgadino, Gualdani and Taskovic \cite{Caja-LopezDelgadinoGualdani2025} obtained the exponential decay of the $L^2$ norm relative to the Maxwellian equilibrium, and proved that this norm is guaranteed to be non-increasing after some time, similarly to Theorem~\ref{thm:main} of the present work. The idea that quantities would become more easily monotonous after an initial thermalization period originated in their work. This thermalization period is easily quantified in the Maxwellian setting thanks to an explicit expression for the directional temperatures.

Of course, the Landau-Maxwell remains a simplified model and we do not claim that the results presented here automatically transfer to more general potentials. Although we are inclined to believe that Theorem~\ref{thm:main} extends to short times, and maybe to more general Landau models, we do not believe our method of proof to be enough to cover those cases without additional new ideas.

\subsection{Notation and layout}
We gather here some notation used in the following sections. The dot product between matrices is
 $$A:B = \Tr(A^TB)=\sum_{ i,j=1}^d A_{ij} B_{ij}.$$
We will mostly apply it with symmetric matrices. The tensor product $v\otimes w$ is the matrix with coefficients $(v\otimes w)_{ij}=v_iw_j$. For a field $S=S(v)$ of symmetric $d\times d$ matrices, we define the generalized Fisher information
\begin{equation}
\label{eq:defI_S}
    I_S (f)=\int_{\mathbb{R}^{d}} f(v) S(v):\left[ \nabla\log f(v)\right]^{\otimes 2}dv = \sum_{i,j=1}^d \int_{\mathbb{R}^{d}} f(v) S_{ij}(v) \partial_i\log f(v) \partial_j\log f(v)dv.
\end{equation}
We will also use the short-hand
$I_i(f) = I_{e_i\otimes e_i}(f)=\int_{\mathbb{R}^{d}} f(\partial_i \log f)^2$,
where $(e_i)_{1\leq i \leq d}$ is the canonical basis of $\mathbb{R}^d$.

For a functional $\mathcal{J}$, we denote its Gâteaux derivative at point $F$ in the direction $H$ by
$$\langle \mathcal{J}(F),H \rangle:=\left.\frac{d}{dt}\mathcal{J}(F_t)\right\vert_{t=0}$$
where $(F_t)$ is a curve such that $F_0=F$ and $\left.\frac{d}{dt}F_t\right\vert_{t=0}=H$.

Finally, the commutator of two operators $A$ and $B$ is $[A,B]:=AB-BA$.

The layout of the remaining sections is as follows: in Section~\ref{sec:Properties of the Landau-Maxwell equation}, we gather rather classical properties and simplifications of the Maxwellian setting that will be used in the sequel. In Section~\ref{sec:A formula for the derivative of the entropy production}, we derive an exact formula for the time-derivative of the entropy production. In Section~\ref{sec:Analysis of the formula and conclusion}, we prove the main results by studying the formula when the directional temperatures are not too unbalanced.

\section*{Acknoledgements}

The author is thankful to Luis Silvestre for sharing with us many insights into kinetic theory, and for pointing out to us a mistake in a previous version of this manuscript, the correction of which lead to the present more direct proof. The author thanks his PhD advisors Cyril Imbert and Clément Mouhot for their advice and discussions around McKean's conjectures. The author finally wishes to thank the organizers of the Kinetic Theory: Novel Statistical, Stochastic and Analytical Methods program during which part of this work was conducted.

\section{Properties of the Landau-Maxwell equation}
\label{sec:Properties of the Landau-Maxwell equation}
We present in this section more or less classical simplifications that occur when considering the Landau equation with Maxwell molecules $\alpha=1$. For the remainder of this note, we fix $f_0$ an initial condition satisfying the normalization \eqref{eq:normalization}, and
consider $f=(f_t)_{t\geq 0}$ the solution of the Landau-Maxwell equation starting from $f_0$. We define the temperature tensor $T$ and a closely related matrix $\chi$ by
\begin{align}
    \label{eq:defT}
        T(t) := \int_{\mathbb{R}^{d}} (v\otimes v)f_t(v) dv, &&\chi(t) := \frac{1}{d-1}\left(\Id-T(t)\right).
\end{align}
The matrix $\chi$ measures the anisotropy of the solution at time $t$, in terms of departure of the temperatures from $1$: for radial solutions, $\chi=0$. Note that $0\leq T\leq d\Id$ by the normalization~\eqref{eq:normalization}, so $-\Id \leq \chi \leq (d-1)^{-1} \Id$.

The simplifications of the Maxwellian setting essentially arise as consequences of the following explicit formula \cite{Villani1998b}: 

\begin{lemma}
    \label{lem:formulaA}
     For all $t\geq 0$, all $v\in\mathbb{R}^d$,
    \begin{equation*}
        \int_{\mathbb{R}^{d}} a(v-w) f_t(w)dw = a(v)+\Id + \chi(t).
    \end{equation*}
\end{lemma}
\begin{proof}
    Observe that
    $$a(v-w)=a(v)+a(w) +\frac{1}{d-1}\left[- 2v\cdot w+v\otimes w + w\otimes v\right].$$
    Multiplying by $f_t(w)$ and integrating in $w$, the bracket vanishes because of the momentum conservation.
    Hence, for any fixed $v$,
    \begin{align*}
        \int_{\mathbb{R}^{d}} a(v-w) f_t(w)dw&=\int_{\mathbb{R}^{d}} \left(a(v)+a(w)\right)f_t(w)dw\\
        &=a(v)+\frac{1}{d-1}\int_{\mathbb{R}^{d}} \left( \vert w\vert^2\Id - w\otimes w\right)f_t(w)dw\\
        &=a(v) +\frac{1}{d-1}\left(d\Id-T(t)\right)
    \end{align*}
    thanks to the energy conservation for the second term. Plugging in the expression of $\chi$ concludes.
\end{proof}
We collect the useful consequences of Lemma~\ref{lem:formulaA} in the remainder of this section. The most important is the following linear reformulation of the Landau equation:

\begin{lemma}
\label{lem:linearlandau}
For $\alpha=1$, the Landau equation~\eqref{eq:landau} rewrites:
\begin{equation}
    \label{eq:landaumax}
    \partial_t f_t = \nabla \cdot\left( (a + \Id + \chi) \nabla f_t \right) + \nabla \cdot(vf_t).
\end{equation}
\begin{proof}
    Since $(\nabla_v - \nabla_w)f(v)f(w) = f(w) \nabla_v f (v) -f(v) \nabla_w f(w)$, equation~\eqref{eq:landau} is the sum of two terms:
    \begin{align*}
        \partial_t f_t(v) = \nabla_v \! \cdot\!\left(  \left(\int_{\mathbb{R}^{d}} a(v-w) f(w)dw\right)\nabla_v f (v)\right)-\nabla_v \! \cdot\! \left( \left(\int_{\mathbb{R}^{d}} a(v-w)\nabla_w f(w) dw\right) f(v)\right). 
    \end{align*}
    By Lemma~\ref{lem:formulaA}, the first term is $\nabla \cdot\left( (a + \Id + \chi) \nabla f_t \right)$. Integrating by parts the second term and using that $\nabla_w\cdot (a(v-w))= v-w$ yields the result.
\end{proof}
\end{lemma}
To complete the above equation, we have the following explicit formula of the temperature tensor and of $\chi$:
\begin{align}
\label{eq:T}
    T(t) = e^{-4dt/(d-1)} (T_0-\Id)+\Id && \chi(t)=e^{-4dt/(d-1)}\chi_0.
\end{align}
We will write $T_{\max}$ for $\max_{1\leq i\leq d} T_i$. The proof of \eqref{eq:T} is a direct computation of $\frac{d}{dt}T_i$ and $\frac{d}{dt}\int v_iv_jf$. A recent reference is \cite[Proposition 2]{Caja-LopezDelgadinoGualdani2025}. The directional temperatures exponentially converging to one correspond to $f_t$ becoming closer to a radial function.

The entropy production can also be simplified, taking the following form:
\begin{align}
\label{eq:formulaD}
    D(f)= I_{a+\Id+\chi}(f) - d.
\end{align}
\begin{proof}[Proof of formula \eqref{eq:formulaD}]
    We differentiate $H(f)=\int f \log f$ along the flow of equation \eqref{eq:landaumax}:
    $$-D(f)=\int \nabla \cdot \left(\left[ a(v)+\Id +\chi \right] \nabla f\right)\log f +\int \nabla\cdot (vf)\log f.  $$
    Integrating the first term by parts yields $-I_{a+\Id+\chi}(f)$. Integrating the second one twice leaves us with
    $\int \nabla\cdot (vf)\log f = -\int v\cdot \nabla f=d\int f=d.$
\end{proof}

The final tool we will use is the decomposition of the matrix field $a$ into vector fields from~\cite{GuillenSilvestre2023}. We introduce, for $1\leq i<j\leq d$, the smooth and divergence free vector fields $b_{ij}$ and their associated differential operator $\mathcal{B}_{ij}$
\begin{align*}
    b_{ij}(v)=(d-1)^{-1/2}\left[v_i e_j - v_j e_i\right] && \mathcal{B}_{ij}f = b_{ij}(v)\cdot \nabla f.
\end{align*}
These vector fields allow for the decomposition
\begin{equation}
    \label{eq:decomp_a}
    a(v)=\sum_{1\leq i<j \leq d}b_{ij}(v)\otimes b_{ij}(v).
\end{equation}
Using that $b_{ij}$ is divergence-free, we can also write $\nabla \cdot(a\nabla f)=\sum_{1\leq i<j \leq d}\mathcal{B}_{ij} \mathcal{B}_{ij}f $.

The flow generated by $b_{ij}$ is a rotation of the plane spanned by $e_i$ and $e_j$. Since the Euclidean Laplacian commutes with rotations, it directly implies that $[\mathcal{B}_{ij},\Delta]=0$ for all $i,j$. However, we stress that the commutators of the vector fields $[\mathcal{B}_{ij},\partial_k]$ is not always zero. Finally, it can be helpful to keep in mind that the operator $(d-1)\nabla \cdot (a\nabla f)$ is exactly the spherical Laplacian of $f$. More precisely, in polar coordinates $v=r\sigma$ with $(r,\sigma)\in \mathbb{R}_+\times \mathbb{S}^{d-1}$, we have $\Delta  f=\frac{1}{r^{d-1}}\partial_r (r^{d-1} \partial_r f) +\frac{1}{r^2}(d-1)\nabla \cdot (a\nabla f)$. This interpretation provides a computation-free proof of $[\mathcal{B}_{ij},\nabla  \cdot a \nabla ]=0$ since the spherical Laplacian commutes with rotations.

\begin{remark}
\label{rem:radialcase}
When $f$ is radial, $a\nabla f=0$ and $\chi=0$ so the Landau equation becomes the Fokker-Planck equation $\partial_t f_t =\Delta f + \nabla\cdot(vf),$ and the entropy production is simply $D=I_{\Id}-d$. It is well-known that the Fisher information is monotone along the flow of this equation. Hence any counter-example to monotonicity could not be radial, making it fundamentally different from the one of \cite{Silvestre2026}.
\end{remark}
\section{Derivative of the entropy production}
\label{sec:A formula for the derivative of the entropy production}
In this section, we compute the time derivative of $D$. For $S$ a symmetric matrix field, we introduce the short-hand notation 
\begin{align*}
    \Delta_Sf:=\nabla\cdot(S\nabla f),
\end{align*}
so that the Landau equation takes the following guise: $\partial_t f = \Delta_{a+\Id+\chi}+\nabla\cdot(vf)$. Recall that $D=I_{a+\Id +\chi}-d$, which depends on time through $\chi$, so we have to differentiate in time the functional $D$ itself:
\begin{equation}
    \label{eq:derivD0}
    \frac{d}{dt}(D(f))=\frac{dD}{dt}(f) + \jap{D'(f),\partial_t f}=-\frac{4d}{d-1}I_\chi(f) + \jap{I_{a+\Id+\chi}'(f),\partial_t f}.
\end{equation}
Here we have used \eqref{eq:T} to write that $\frac{d}{dt}I_\chi= -\frac{4d}{d-1}I_\chi$. We now have to compute the directional derivative above, beginning with the drift part:
\begin{lemma} 
\label{lem:drift}
It holds that
$$ \jap{I_{a+\Id +\chi}'(f),\nabla\cdot(vf)}= 2I_{\Id + \chi}(f).$$
\end{lemma}
\begin{proof}
    Let $A(v)=a(v)+Id+\chi$. The solution of $\partial_t f_t = \nabla\cdot(vf_t)$ with initial data $f$ is $f_t(v) = e^{dt}f(e^tv)$, which satisfies
    $\nabla \log f_t (v)= e^t(\nabla \log f) (e^t v)$. Hence, $$I_{A}(f_t)= e^{2t}\int e^{dt} f(e^{t} v)A(v): \left[(\nabla \log f)(e^{t}v)\right]^{\otimes 2} dv =e^{2t}\int  f( v)A(e^{-t}v) \left[\nabla \log f (v)\right]^{\otimes 2} dv$$
    by changing variables. We now plug in that $A(e^{-t}v)=e^{-2t}a(v)+\Id+\chi$, getting that $I_A(f_t)=I_a(f)+e^{2t}I_{\Id +\chi}(f)$. We get the result by taking the time derivative of this expression at $t=0$.
\end{proof}

It remains to compute $\jap{I_{a+\Id+\chi}'(f),\Delta_{a+\Id+\chi} f}$, which is the sum of nine possible pairings. On one hand, since $\Delta_a$ and $\Delta=\Delta_{\Id}$ commute, it will be easy to compute the derivative of $I_{a}$ in the direction $\Delta_{\Id}$ and conversely, and see that it is non-negative. On the other hand, $\Delta_\chi$ also commutes with $\Delta$ since $\chi$ is constant in $v$, but not with $\Delta_a$: the interaction of $a$ and $\chi$ leads to the most complicated terms.

To proceed with the computation, we provide two abstract lemmas. They turn out most useful when the operators at play commute, but we state them in full generality.
\begin{lemma}
\label{lem:abs1}
    Consider the matrix field $S(v)=b(v)\otimes b(v)$ for some smooth vector field $b$ and the associated derivation operator $\mathcal{B}=b(v)\cdot \nabla$. For any self-adjoint operator $L$, define the associated \textit{carré-du-champ}
    $$2\Gamma^L(g)=Lg^2-2(Lg)g.$$
    Then,
    $$\jap{I_S'(f),Lf}=-\int 2f\Gamma^L(\mathcal{B} \log f)  dv+2\int(\mathcal{B}\log f)[\mathcal{B},L]f.$$
\end{lemma}
\begin{proof}
    Since $I_S(f) = \int ( \mathcal{B} f)^2/f \ dv$, a direct computation yields
    \begin{align*}
        \jap{I_S'(f),Lf} &= -\int (Lf)(\mathcal{B} \log f)^2 +\int  2 (\mathcal{B}\log f )(\mathcal{B} Lf)\\
       & = -\int f L( \mathcal{B} \log f)^2 +\int  2 (\mathcal{B} \log f )(L\mathcal{B} f) + 2\int( \mathcal{B}\log f)[\mathcal{B},L]f
    \end{align*}
    by self-adjointness of $L$ for the first integrand. We keep the third term untouched, and observe that the middle term rewrites
    $$\int 2 (\mathcal{B} \log f )(L\mathcal{B} f)=\int  2( L\mathcal{B}\log f )(\mathcal{B} f)=\int  2 f(L\mathcal{B} \log f )(\mathcal{B}  \log f).$$
    Combining it with the first term yields the result.
\end{proof}
When the operator $L$ can also be expressed with vector fields, its \textit{carré-du-champ} is non-negative with a simple expression:
\begin{lemma}
\label{lem:abs2}
    Consider a smooth divergence-free vector field $\tilde{b}$ and define the second order operator $L=(\tilde{b}\cdot \nabla)( \tilde{b}\cdot \nabla).$ Then $L$ is formally self-adjoint and its carré-du-champ is given by $\Gamma^L(g) = ( \tilde{b}\cdot \nabla g)^2$.
\end{lemma}
\begin{proof}
    The self-adjointness is shown by integrating by parts the operator $\tilde{b}\cdot \nabla$ twice, which is made possible thanks to the divergence-free condition. For the \textit{carré-du-champ}, the direct computation
    $$Lg^2=\tilde{b}\cdot \nabla(2g\tilde{b}\cdot \nabla g)=2g(Lg)+2(\tilde{b}\cdot \nabla g)^2$$
    gives the result.
\end{proof}
Applying these Lemmas allows us to compute all the remaining terms but one. Recall, from the end of Section~\ref{sec:Properties of the Landau-Maxwell equation}, the commutation $[\mathcal{B}_{ij},\Delta_{\Id}]=0$ for all $i<j$, which entails $[\Delta_a,\Delta_{\Id}]=0$. Also, $[\mathcal{B}_{ij},\Delta_{a}]=0$. In general, the commutator of $\mathcal{B}_{ij}$ with partial derivatives is given by:
\begin{align}
\label{eq:commutator}
    [\partial_i,\mathcal{B}_{kl}]=(d-1)^{-1/2}\left(\delta_{ik}\partial_l-\delta_{il}\partial_k\right)=(d-1)^{-1/2}\begin{cases}
        \partial_l & \text{if }i=k \\
        -\partial_k & \text{if }i=l \\
        0 & \text{else.}
    \end{cases}
\end{align}
We are now poised to show:
\begin{lemma} 
\label{lem:niceterms}
It holds that
\begin{align}
\label{eq:I'a,a}
    \jap{I'_{a}(f),\Delta_a f)}&=-2\sum_{i<j,k<l}\int f \vert \mathcal{B}_{ij}\mathcal{B}_{kl}\log f\vert^2,\\
\label{eq:I'a,Id}
    \jap{I'_{a}(f),\Delta f}=\jap{I'_{\Id}(f),\Delta_a f}&=-2\sum_{k<l}\int f \vert \nabla\mathcal{B}_{kl} \log f\vert^2,\\
\label{eq:I'id+chi,id+chi}
    \jap{I'_{\Id+\chi}(f),\Delta_{\Id+\chi} f}&=-2\sum_{i,j}(1+\chi_i)(1+\chi_j)\int f \vert \partial_{ij}^2 \log f\vert^2,\\
\label{eq:I'a,chi}
    \jap{I'_a(f),\Delta_\chi f}&=-\sum_i\!\sum_{k<l}\chi_i\bigg[2\int f(\partial_i\mathcal{B}_{kl} \log f)^2\\
    \nonumber
    &\qquad +4\int f(\partial_i\mathcal{B}_{kl} \log f)([\mathcal{B}_{kl},\partial_i]\log f) \bigg].
\end{align}
All the indices in the sums range from $1$ to $d$.
\end{lemma}

\begin{proof}
\textit{Formula~\eqref{eq:I'a,a}.} We apply the decomposition \eqref{eq:decomp_a} of $a$ to get $I_a=\sum_{i<j} I_{b_{ij}\otimes b_{ij}}$. To each term, we apply Lemma~\ref{lem:abs1} with $\mathcal{B}=\mathcal{B}_{ij}$ and $L=\Delta_a$. The commutator vanishes, so
\begin{align*}
    \jap{I'_{a}(f),\Delta_a f} =-\sum_{i<j}\int 2f\Gamma^{\Delta_a}(\mathcal{B}_{ij} \log f ).
\end{align*}
But using that $\Delta_a=\sum_{k<l}\mathcal{B}_{kl}\mathcal{B}_{kl}$, Lemma~\ref{lem:abs2} then tells us that $\Gamma^{\Delta_a}(g ) =\sum_{k<l}(\mathcal{B}_{kl}g)^2 $. This yields formula~\eqref{eq:I'a,a}.

\textit{Formula~\eqref{eq:I'a,Id}.}
We first show the equality between the two derivatives. Recall that $I_{\Id} (f)$ is minus the derivative of $H$ along the Euclidean heat flow: $I_{\Id} (f) = -\jap{H'(f),\Delta f}$. Similarly, $I_a(f)=-\jap{H'(f), \Delta_a f}$. Since $\Delta_a$ and $\Delta$ commute, differentiating $H$ along $\Delta$ and then along $\Delta_a$ is the same as differentiating first along $\Delta_a$ and then along $\Delta$. This yields $\jap{I'_{a}(f),\Delta f}=\jap{I'_{\Id}(f),\Delta_a f}$.

We now show the second half of the formula. Applying again the decomposition \eqref{eq:decomp_a} of $a$, Lemma~\ref{lem:abs1} with $\mathcal{B}=\mathcal{B}_{ij}$ and $L=\Delta$ gives
$$\jap{I'_{a}(f),\Delta f} = -\sum_{i<j}\int 2f \Gamma^\Delta(\mathcal{B}_{ij} \log f ),$$
because the commutator vanishes again. Since $\Gamma^\Delta(g)=\vert \nabla g\vert^2$ (by Lemma~\ref{lem:abs2} for instance) we get Formula~\eqref{eq:I'a,Id}.

\textit{Formula~\eqref{eq:I'id+chi,id+chi}.} 
We decompose $I_{Id+\chi}=\sum_i (1+\chi_i) I_{e_i\otimes e_i}$ and $\Delta_{\Id+\chi} f =\sum_j (1+\chi_j)\partial_{jj}^2f$ on the canonical basis to get
\begin{align*}
    \jap{I'_{\Id+\chi}(f),\Delta_{\Id+\chi} f}&=\sum_{i,j}(1+\chi_i)(1+\chi_j)\jap{I'_{e_i\otimes e_i}(f),\partial^2_{jj} f}.
\end{align*}
Of course $\partial_i$ and $\partial_{jj}^2$ commute, and Lemmas~\ref{lem:abs1} and~\ref{lem:abs2} prove the claimed formula.

\textit{Formula~\eqref{eq:I'a,chi}.} We use the decomposition~\eqref{eq:decomp_a} and $\Delta_\chi=\sum_i \chi_i\partial^2_i$ to write
    $$\jap{I'_a(f),\Delta_\chi f} = \sum_i \sum_{k<l} \chi_i \jap{I'_{b_{kl}\otimes b_{kl}}(f),\partial^2_{ii} f}.$$
This time there is commutator terms. Observe from~\eqref{eq:commutator} that $[\mathcal{B}_{kl},\partial_i]$ commutes with partial derivatives (because it is one itself), so $[\mathcal{B}_{kl},\partial^2_{ii}]=\partial_i[\mathcal{B}_{kl},\partial_{i}]+[\mathcal{B}_{kl},\partial_{i}]\partial_i=2\partial_i[\mathcal{B}_{kl},\partial_{i}]$. Applying Lemmas~\ref{lem:abs1} and~\ref{lem:abs2},
$$\jap{I'_{b_{kl}\otimes b_{kl}}(f),\partial^2_{ii} f}=-2\int f(\partial_i \mathcal{B}_{kl} \log f)^2 +4 \int (\mathcal{B}_{kl} \log f)\partial_i[\mathcal{B}_{kl},\partial_{i}] f.$$
Integrating by part the $\partial_i$ in the last term yields the result.
\end{proof}
\begin{remark}
    If we try to compute $\jap{I'_{\Id}(f),\Delta_a f}$ directly rather than writing it as $\jap{I'_a(f),\Delta f}$, the commutator $[\partial_i,\Delta_a]$ is not zero, so we obtain the expression:
    $$\jap{I'_{\Id}(f),\Delta_a f}=-2\sum_{k<l}\int f \vert \mathcal{B}_{kl} \nabla \log f\vert^2+2\sum_{k<l}\int(\nabla\log f)[\nabla,\Delta_a]f,$$
    which is less easy to work with. In particular, notice that the left-hand side vanishes on radial functions because then $\Delta_a f=0$, which is not clear on the right-hand side.
\end{remark}
The only remaining term is $\jap{I'_\chi(f),\Delta_a f}$. If we apply Lemmas~\ref{lem:abs1} and~\ref{lem:abs2} to compute it, it features the square $(\mathcal{B}_{kl} \partial_i \log f)^2$ and a complicated commutator term. We would rather have the opposite order of operators, and a cleaner error term: for this we use this third abstract Lemma:
\begin{lemma}
\label{lem:abs3}
    Let $\mathcal{B}=b\cdot\nabla$, $\mathcal{L}=\tilde{b}\cdot\nabla$ be the two derivation operators associated to divergence free vector fields $b$ and $\tilde{b}$. Then,
    \begin{align*}
        \jap{I'_{b\otimes b}(f),\mathcal{L}\mathcal{L}f}=-2\int f (\mathcal{B}\mathcal{L}(\log f))^2 + 2\int f ([\mathcal{B},\mathcal{L}]\log f)^2 +2\int f(\mathcal{B}\log f)([[\mathcal{B},\mathcal{L}],\mathcal{L}]\log f)
    \end{align*}
\end{lemma}
\begin{proof}
    We start from $I_{b\otimes b}(f)=\int (\mathcal{B}f)^2/f$ and differentiate along $\mathcal{L}$. The goal is to bring $\mathcal{L}$ to the front using commutators to get a full derivative with respect to $\mathcal{L}$. We proceed as follows:
    \begin{align*}
        \jap{I'_{b\otimes b}(f),\mathcal{L}f}&= \int \left[ 2(\mathcal{B} f) (\mathcal{B} \mathcal{L}f)\frac{1}{f}+(\mathcal{B}f)^2\mathcal{L}\left(\frac{1}{f}\right)\right]\\
        &=\int \left[ 2(\mathcal{B} f) ([\mathcal{B} ,\mathcal{L}]f)\frac{1}{f}+2(\mathcal{B} f) (\mathcal{L}\mathcal{B}f)\frac{1}{f}+(\mathcal{B}f)^2\mathcal{L}\left(\frac{1}{f}\right)\right]\\
        &=2\int f (\mathcal{B} \log f) ([\mathcal{B} ,\mathcal{L}]\log f) + \int\mathcal{L} \left[\frac{(\mathcal{B}f)^2}{f}\right].
    \end{align*}
    The second term vanishes by the divergence-free condition. We differentiate a second time, applying the same method to obtain:
    \begin{align*}
        \jap{I'_{b\otimes b}(f),\mathcal{L}\mathcal{L}f}+\jap{I''_{b\otimes b}(f),\mathcal{L}f,\mathcal{L} f}=2\int f  ([\mathcal{B} ,\mathcal{L}]\log f)^2+2\int f (\mathcal{B} \log f) ([[\mathcal{B} ,\mathcal{L}]\mathcal{L}]\log f).
    \end{align*}
    We check that the second derivative of $I_{b\otimes b}$ is given by $\jap{I''_{b\otimes b}(f),h,h}=2\int f (\mathcal{B}(h/f))^2,$ which is a straightforward computation.
\end{proof}
Applying this lemma yields the final derivative:
\begin{lemma}
\label{lem:badterms}
    It holds that:
    \begin{align}
    \label{eq:I'chi,a}
        \jap{I'_{\chi}(f),\Delta_a f}&=-2\sum_i\!\sum_{k<l}\chi_i \int f(\partial_i \mathcal{B}_{kl}\log f)^2 -\frac{2d}{d-1}I_\chi(f).
    \end{align}
\end{lemma}
\begin{proof}
    We write the decomposition
    $$ \jap{I'_\chi(f),\Delta_a f} = \sum_i \sum_{k<l} \chi_i \jap{I'_{e_i\otimes e_i}(f),\mathcal{B}_{kl}\mathcal{B}_{kl}f}$$
    and apply Lemma~\ref{lem:abs3}:
    \begin{align*}
         \jap{I'_\chi(f),\Delta_a f} = \sum_i \sum_{k<l} \chi_i\bigg[-2\int f(\partial_i \mathcal{B}_{kl}&\log f)^2+2\int f ([\partial_i,\mathcal{B}_{kl}]\log f)^2\\
         &+2\int f(\partial_i\log f)([[\partial_i,\mathcal{B}_{kl}],\mathcal{B}_{kl}]\log f)\bigg]
    \end{align*}
    The first term we keep as it is. The second one, using the expression \eqref{eq:commutator} for $[\partial_i,\mathcal{B}_{kl}]$, is
    $$\sum_i\sum_{k<l}\chi_i([\partial_i,\mathcal{B}_{kl}]\log f)^2 = \frac{1}{d-1}\sum_{k<l}\chi_k (\partial_l\log f)^2+\chi_l (\partial_k\log f)^2= \frac{1}{d-1}\sum_{k\neq l}\chi_k (\partial_l\log f)^2$$
    using that the terms under the sum are symmetric in $k,l$.
    But since $\sum_k \chi_k = \text{Tr}(\chi)= 0$, we have for any fixed $l$, $\sum_{k\  \text{s.t.}\ k\neq l} \chi_k=-\chi_l$, leading to
    $$ \sum_i\sum_{k<l}2\chi_i\int f ([\mathcal{B}_{kl},\partial_i ]\log f)^2=-\frac{2}{(d-1)}\sum_{l}\chi_l \int f(\partial_l\log f)^2=-\frac{2}{d-1} I_\chi(f).$$
    
    For the last term, we use that
    $$[[\partial_i,\mathcal{B}_{kl}],\mathcal{B}_{kl}]=\frac{1}{d-1}\begin{cases}
        -\partial_k & \text{if }i=k \\
        -\partial_l & \text{if }i=l \\
        0 & \text{else.}\end{cases} .$$
    This implies
    \begin{align*}
        \sum_i\sum_{k<l} \chi_i(\partial_i \log f) ([[\partial_i,\mathcal{B}_{kl}],\mathcal{B}_{kl}]\log f)&=-\frac{1}{d-1}\sum_{k<l}\chi_k(\partial_k \log f)^2+\chi_l(\partial_l \log f)^2 \\
        &= -\frac{1}{d-1}\sum_{k\neq l}\chi_k(\partial_k \log f)^2=-\sum_k \chi_k(\partial_k \log f)^2,
    \end{align*}
    by symmetry. Thus
    $$\sum_i \sum_{k<l}2\chi\int f(\partial_i\log f)([[\partial_i,\mathcal{B}_{kl}],\mathcal{B}_{kl}]\log f)=-2I_\chi(f),$$
    hence the result.
\end{proof}
\begin{remark}
    Using Lemma~\ref{lem:abs3} and proceeding similarly, we can also show the formula
    $$\jap{I'_{a}(f),\Delta_\chi f}=-2\sum_i\!\sum_{k<l}\chi_i \int f(\mathcal{B}_{kl}\partial_i\log f)^2 -\frac{2}{d-1}I_\chi(f),$$
    which is somewhat cleaner than~\eqref{eq:I'a,chi}. However, it features the order of operators $\mathcal{B}_{kl}\partial_i$ instead of $\partial_i\mathcal{B}_{kl}$.
\end{remark}

Gathering everything we have shown, we can obtain an organized formula for $\frac{d}{dt}D(f)$. Starting from \eqref{eq:derivD0},
\begin{align*}
    \frac{d}{dt}D(f)&=-\frac{4d}{d-1}I_\chi(f)+\jap{I_{a+\Id+\chi}'(f),\nabla \cdot (vf)}+\jap{I_{a+\Id+\chi}'(f),\Delta_{a+\Id+\chi}  f}.
\end{align*}
Using Lemma~\ref{lem:drift} for the drift term, 
\begin{align}
\label{eq:derivD1}
    \frac{d}{dt}D(f)
    &=-\frac{2(d+1)}{d-1}I_\chi(f)+2I_{\Id}(f)+\jap{I_{a+\Id+\chi}'(f),\Delta_{a+\Id+\chi}  f}.
\end{align}
We expand the last term:
\begin{align*}
   \jap{I_{a+\Id+\chi}'(f),\Delta_{a+\Id+\chi}  f}  &=\jap{I_{a}'(f),\Delta_a  f}\\
    &+\jap{I_{\Id}'(f),\Delta_a  f}+\jap{I_{a}'(f),\Delta  f}\\
    &+\jap{I'_{\Id+\chi}(f),\Delta_{\Id+\chi} f}\\
    &+\jap{I'_{a}(f),\Delta_{\chi}f}+\jap{I'_{\chi}(f),\Delta_{a} f}.
\end{align*}
The two terms on the second line are the same by formula~\eqref{eq:I'a,Id}, and can be combined with parts of the two terms of the last line, given by~\eqref{eq:I'a,chi} and~\eqref{eq:I'chi,a}. The first and third lines are given by formulas~\eqref{eq:I'a,a} and \eqref{eq:I'id+chi,id+chi} respectively. We obtain:
\begin{align*}
   \jap{I_{a+\Id+\chi}'(f),\Delta_{a+\Id+\chi}  f}=&-2\sum_{i<j,k<l}\int f \vert \mathcal{B}_{ij}\mathcal{B}_{kl}\log f\vert^2\\
   &-4\sum_i\sum_{k<l}(1+\chi_i)\int f (\partial_i\mathcal{B}_{kl} \log f)^2\\
   &-2\sum_{i,j}(1+\chi_i)(1+\chi_j)\int f \vert \partial_{ij}^2 \log f\vert^2\\
   &-4\sum_i\sum_{k<l} \chi_i\int f(\partial_i\mathcal{B}_{kl} \log f)([\mathcal{B}_{kl},\partial_i]\log f)-\frac{2d}{d-1}I_\chi(f).
\end{align*}
Introducing the notation
$$K^{\mathcal{B}\mathcal{L}}(f)=2\int f (\mathcal{B}\mathcal{L} \log f)^2$$
and gathering everything, we have proven:
\begin{prop}[Non-relative formula for $\frac{d}{dt}D(f)$]
The derivative of the entropy production is given by
\begin{align}
\label{eq:derivD2}
    \frac{d}{dt}&D(f)
    =-\frac{2(2d+1)}{d-1}I_\chi(f)+2I_{\Id}(f)+\mathcal{E}\\
    &-\sum_{i<j,k<l}K^{ \mathcal{B}_{ij}\mathcal{B}_{kl}}(f)-2\sum_i\sum_{k<l}(1+\chi_i)K^{\partial_i\mathcal{B}_{kl} }(f)-\sum_{i,j}(1+\chi_i)(1+\chi_j)K^{\partial_i\partial_j}(f), \nonumber
\end{align}
with the errror
\begin{equation}
    \label{eq:error}
    \mathcal{E}=-4\sum_i\sum_{k<l}\chi_i \int f(\partial_i\mathcal{B}_{kl} \log f)([\mathcal{B}_{kl},\partial_i]\log f).
\end{equation}
\end{prop}

Since $\chi \geq -1$, the second line of \eqref{eq:derivD2} is non-positive. However, even in the radial case, this expression of $D$ is still not obviously non-positive. The sign becomes most obvious when rewriting the terms relatively to the Maxwellian $m=(2\pi)^{d/2}e^{-\vert v\vert^2/2}$.
We hence define:
\begin{align*}
    I_S(f\vert m) := \int f (v)S(v):\left[\nabla \log\left(\frac{f}{m}\right)\right]^{\otimes 2},&&K^{\mathcal{B}\mathcal{L}}(f\vert m) :=2\int f\left(\mathcal{B}\mathcal{L}\log\left(\frac{f}{m}\right)\right)^2.
\end{align*}
We also use the short-hand $I_i(f\vert m)=I_{e_i\otimes e_i}(f\vert m)$. The relation between relative and non-relative quantities is as follows:
\begin{lemma}
\label{lem:relat}
    It holds that
    \begin{align*}
        I_{i}(f) &=  I_{i}(f\vert m)-T_i +2, \\
        K^{\partial_i\partial_j}(f)
        &= K^{\partial_i\partial_j}(f\vert m)+4\delta_{ij}I_{i}(f\vert m)+2\delta_{ij}(3-2T_i),
    \end{align*}
leading to the formulas
\begin{align*}
    I_{\Id}(f)&=I_{\Id}(f\vert m)+d, && I_\chi(f)=I_\chi(f\vert m)+(d-1)\Vert \chi\Vert^2,
\end{align*}
and
\begin{align*}
    \sum_{i,j}(1+\chi_i)(1+\chi_j)K^{\partial_i\partial_j}(f)=&\sum_{i,j}(1+\chi_i)(1+\chi_j)K^{\partial_i\partial_j}(f\vert m)+4I_{(\Id + \chi)^2}(f\vert m)\\
    &+2d+(8d-6)\Vert \chi\Vert^2+4(d-1)\sum_i\chi_i^3,
\end{align*}
where $\Vert \chi\Vert^2=\sum_i \chi_i^2$.
\end{lemma}
\begin{proof}
    For the first formula, we use that $\partial_i\log m=-v_i$ to write:
    \begin{align*}
        I_{i}(f\vert m)=\int f\left(\partial_i\log\left(f\right)+v_i\right)^2 &=I_{i}(f)+2\int( \partial_if) v_i+\int f v_i^2=I_{i}(f)-2+T_i,
    \end{align*}
    by integrating by parts the middle term.
    Similarly for the second formula, $\partial_{ij}^2\log m=-\delta_{ij}$, so
    \begin{align*}
         2\int f(\partial_{ij}^2\log (f/m) )^2 &=  2\int f(\partial_{ij}^2\log f )^2+4\delta_{ij}\int f(\partial_{ii}^2\log f)+2\delta_{ij}\int f\\
         &=2\int f(\partial_{ij}^2\log f )^2-4\delta_{ij}\int f(\partial_{i}\log f)^2+2\delta_{ij}.
    \end{align*}
    We then use the previous expression to rewrite $\int f(\partial_{i}\log f)^2=I_i(f)$ using relative information and obtain the second formula.

    Summing the first formula over $i$ and recalling $\sum T_i=d$ yields the third one. Also, $$I_\chi(f)=\sum_i \chi_i I_{i}(f\vert m)-\chi_i T_i+2\chi_i= I_\chi(f\vert m)+(d-1)\sum_i\chi_i^2$$
using that the $\chi_i$ sum to $0$, and that $T_i=1-(d-1)\chi_i$. Finally,
\begin{align*}
    \sum_{i,j}(1+\chi_i)(1+\chi_j)K^{\partial_i\partial_j}(f)=&\sum_{i,j}(1+\chi_i)(1+\chi_j)K^{\partial_i\partial_j}(f\vert m)\\
    &+4\sum_{i}(1+\chi_i)^2I_{i}(f\vert m) + 2\sum_i{(1+\chi_i)^2}(3-2T_i).
\end{align*}
We further compute the constant:
\begin{align*}
    2\sum_i{(1+\chi_i)^2}(3-2T_i)&=2\sum_i{(1+\chi_i)^2}(1+2(d-1)\chi_i)\\
    &=2d+(8d-6)\Vert \chi \Vert^2 +4(d-1)\sum_i\chi_i^3,
\end{align*}
using that the $\chi_i$ sum to $0$.
\end{proof}

Recall that $b_{kl}$ generates a rotation, so for any radial function $g$, $\mathcal{B}_{kl}g=0$. In particular, since $m$ is radial, $K^{ \mathcal{B}_{ij}\mathcal{B}_{kl}}(f)=K^{ \mathcal{B}_{ij}\mathcal{B}_{kl}}(f\vert m)$ and $K^{\partial_i\mathcal{B}_{kl} }(f)=K^{\partial_i\mathcal{B}_{kl} }(f\vert m)$. Plugging this and the results of the previous lemma in, we obtain the relative expression:
\begin{prop}[Relative formula for $\frac{d}{dt}D(f)$]
The derivative of the entropy production is given by
\begin{align}
\label{eq:derivD3}
    \frac{d}{dt}D(f)
    =-&\frac{2(2d+1)}{d-1}I_\chi(f\vert m)+2I_{\Id}(f\vert m)-4I_{(\Id +\chi)^2}(f\vert m)+\mathcal{E}\\
    &-\sum_{i<j,k<l}K^{ \mathcal{B}_{ij}\mathcal{B}_{kl}}(f\vert m)-2\sum_i\sum_{k<l}(1+\chi_i)K^{\partial_i\mathcal{B}_{kl} }(f\vert m)\nonumber\\
    &-\sum_{i,j}(1+\chi_i)(1+\chi_j)K^{\partial_i\partial_j}(f\vert m)\nonumber\\
    &-(12d-4)\Vert \chi\Vert^2-4(d-1)\sum_i\chi_i^3,\nonumber
\end{align}
with the errror $\mathcal{E}$ given by~\eqref{eq:error}.
\end{prop}
Now we retrieve clearly that for $\chi=0$, the entropy production is monotone, as the expression becomes
$$\frac{d}{dt}D(f)=-2I_{\Id}(f\vert m)-\sum_{i,j}K^{\partial_i \partial_j}(f\vert m),$$
which is exactly the dissipation of Fisher information along the Fokker-Planck equation (see Remark~\ref{rem:radialcase}).

\section{Monotonicity and decay for small $\chi$}
\label{sec:Analysis of the formula and conclusion}

Analyzing the formula, we are ready to show monotonicity of $D$ for small enough $\chi$. Remark that when $\chi=0$, we have some room to play in terms of Fisher information, but we still need to be careful of not introducing any positive constant term that we could not control as $\chi\rightarrow 0$ (such a linear term in $\chi$ for instance, as we only have higher powers of $\chi$ in the constant terms of \eqref{eq:derivD3}).
 
First we provide a blunt bound for the error term $\mathcal{E}$.
\begin{lemma}
    \label{lem:errorbound}
    For any $\eta>0$,
    \begin{align*}
        \vert \mathcal{E}\vert=\bigg\vert\sum_{k<l} \sum_i 4\chi_i\int f(\partial_i\mathcal{B}_{kl} \log f)([\mathcal{B}_{kl},\partial_i]&\log f) \bigg\vert\leq \eta\sum_{k<l}\sum_i K^{\partial_i\mathcal{B}_{kl} }(f\vert m)\\       
        &+\frac{2\eta^{-1}}{d-1} \Vert \chi \Vert^2 I_{\Id}(f\vert m)+\frac{2d\eta^{-1}}{d-1} \Vert \chi \Vert^2. 
    \end{align*}
\end{lemma}
\begin{proof}
    We use that $2ab\leq \eta a^2 +\eta^{-1}b^2 $, giving us
    $$4\chi_i\int f(\partial_i\mathcal{B}_{kl} \log f)([\mathcal{B}_{kl},\partial_i]\log f) \leq 2\eta\int f(\partial_i\mathcal{B}_{kl} \log f)^2+2\eta^{-1}\chi_i^2\int f([\mathcal{B}_{kl},\partial_i]\log f)^2$$
    Summing over $k<l$, and over $i$, we obtain using the value of the commutators:
    \begin{align*}
        \sum_{k<l} &\sum_i 4\chi_i\int f(\partial_i\mathcal{B}_{kl} \log f)([\mathcal{B}_{kl},\partial_i]\log f)\\ \leq &\eta \sum_{k<l}\sum_i K^{\partial_i\mathcal{B}_{kl} }(f\vert m) +\frac{2\eta^{-1}}{d-1}\sum_{k<l} \left(\chi_k^2\int f(\partial_l\log f)^2+\chi_l^2\int f(\partial_k\log f)^2\right).
    \end{align*}
    The first  term is what we seek. The last sum rewrites, by symmetry, as:
    $$\frac{2\eta^{-1}}{d-1}\sum_{k\neq l} \chi_k^2\int f(\partial_l\log f)^2\leq \frac{2\eta^{-1}}{d-1}\sum_{k, l} \chi_k^2\int f(\partial_l\log f)^2=\frac{2\eta^{-1}}{d-1} \left(\sum_{k} \chi_k^2\right)I_{\Id}(f).$$
    Writing $I_{\Id}(f)=I_{\Id}(f\vert m)+d$ concludes the proof.
\end{proof}

We are now ready to prove Theorem~\ref{thm:main}.
\begin{proof}[Proof of Theorem~\ref{thm:main}]
We start from the relative form of $\frac{d}{dt}D(f)$ expressed in \eqref{eq:derivD3}. We drop two out of the three non-negative $K$ terms, and plug in the result of the previous lemma for some $\eta>0$ to fix later:
\begin{align}
\label{eq:derivD4}
    \frac{d}{dt}D(f)
    \leq-&\frac{2(2d+1)}{d-1}I_\chi(f\vert m)+2\left[1+\frac{\eta^{-1}}{d-1}\Vert \chi\Vert^2\right]I_{\Id}(f\vert m)-4I_{(\Id +\chi)^2}(f\vert m) \\
    &-2\sum_i\sum_{k<l}(1+\chi_i-\eta/2)K^{\partial_i\mathcal{B}_{kl} }(f\vert m)\nonumber\\
    &+\left[\frac{2d\eta^{-1}}{d-1}-(12d-4)\right]\Vert \chi\Vert^2-4(d-1)\Vert \chi\Vert^3\nonumber.
\end{align}
Since $\vert \chi_i\vert\leq 1$, we have $\vert \sum_i \chi_i^3 \vert \leq \Vert \chi\Vert^2$, and the last line is bounded by
$$\left[\frac{2d\eta^{-1}}{d-1}-(12d-4)\right]\Vert \chi\Vert^2-4(d-1)\Vert \chi\Vert^3 \leq \left[\frac{2d\eta^{-1}}{d-1}-8d\right]\Vert \chi\Vert^2.$$
If we pick $\eta\geq \frac{1}{4(d-1)}$, the right-hand side is non-positive, and the last line of \eqref{eq:derivD4} is non-positive for all times. For the second line, we need to pick $\eta<2$, and it is then non-positive as soon as
\begin{align*}
    \max_i(T_i)-1=(d-1)\max_i(-\chi_i)\leq (d-1)(1-\eta /2)=:C_1(d,\eta)
\end{align*}
The first line is non-positive as soon as the matrix
\begin{align*}
    -\left[\frac{2(2d+1)}{d-1}+
    8\right]\chi +\left[-2+\frac{2\eta^{-1}}{d-1}\Vert \chi\Vert^2\right] \Id -4\chi^2\leq 0.
\end{align*}
This always happens since the left-hand side goes to $-2\Id$ as $\chi\rightarrow 0$. To get a quantitative bound, since $$\Vert\chi\Vert^2=\frac{1}{(d-1)^2}\left(\sum_i T_i^2-d\right)\leq \frac{d}{(d-1)^2}\left((\max_i(T_i))^2-1\right)\leq \frac{d(d+1)}{(d-1)^2}\left((\max_i(T_i))-1\right),$$ we can crudely bound every (diagonal) entry of the left-hand side by:
$$ \left[\frac{2(2d+1)}{d-1}+
    8\right]\frac{\max_i(T_i)-1}{d-1} +\frac{2\eta^{-1}d(d+1)}{(d-1)^3}\left( (\max_i(T_i))-1\right)-2.$$
Multiplying by $(d-1)^2$, we see that this is non-positive as soon as
$$\max_i(T_i)-1\leq (d-1)^2\left((6d-3)+\frac{\eta^{-1}d(d+1)}{d-1}\right)^{-1} =:C_2(d,\eta).$$
We now optimize $\eta\in[\frac{1}{4(d-1)},2)$ to get matching constants $C_1$ and $C_2$.
\begin{itemize}
    \item For $d=2$, we can pick $\eta=\frac{5+\sqrt{133}}{9}\approx 1.83$ so that $C_1(2,\eta)=C_2(2,\eta)=\frac{13-\sqrt{133}}{18}$. 
    \item For $d=3$, we can pick
$\eta=\frac{10+2\sqrt{70}}{15}\approx 1.78$ so that $C_1(3,\eta)=C_2(3,\eta)=\frac{20-2\sqrt{70}}{15}$.
    \item As $d\rightarrow \infty$, $C_1(d,\eta) \sim d(1-\eta/2) $ and $C_2(d,\eta) \sim d/(6+\eta^{-1}) $, so we take $\eta=\frac{9+\sqrt{129}}{12}\approx 1.69$, for which $C_1(d,\eta)\sim C_2(d,\eta)$, and for which one can check that $C_1(d,\eta)\geq C_2(d,\eta)$ always.
\end{itemize}
In any case, the entropy production is non-increasing as soon as $\max_i(T_i)-1\leq C_2(d,\eta)$. If this is not the case initially, the evolution of $T$ given by \eqref{eq:T} tells us that it is guaranteed as soon as
\begin{align*}
    e^{-4dt/(d-1)}(T_{0,\max}-1)\leq C_2(d,\eta), && \textit{i.e.} && t\geq \frac{d-1}{4d}\log\left(\frac{T_{0,\max}-1}{C_2(d,\eta)}\right).
\end{align*}
This proves Theorem~\ref{thm:main} with $C_d=C_2(d,\eta)$.
\end{proof}

We now turn to the:
\begin{proof}[Proof of Theorem~\ref{thm:side}]
    For a function $g$ on $\mathbb{S}^{d-1}$, the $\Gamma^2$ criterion on the sphere \cite{BakryGentilLedoux2014} writes
$$-\frac{1}{2}\jap{\mathcal{I}'(g),\Delta_\sigma g}\geq (d-1)\mathcal{I}(g),$$
where $\Delta_\sigma$ is the spherical Laplacian, and $\mathcal{I}(g) = \int_{\mathbb{S}^{d-1}} \vert \nabla_\sigma \log g\vert^2 g=-\int_{\mathbb{S}^{d-1}} (\Delta_\sigma g)\log g$ is the spherical Fisher information. If we write $f$ in spherical coordinates $f(r,\sigma)$, and recall that $\Delta_\sigma f= (d-1)\Delta_a f$, we get that
$$\int_0^\infty \mathcal{I}(f(r,\cdot))r^{d-1}dr=-(d-1)\int_0^\infty \int_{\mathbb{S}^{d-1}} (\Delta_a f)(\log f) d\sigma r^{d-1}dr=(d-1)I_a(f).$$
Hence, integrating the criterion over $r$, the $(d-1)^2$ factor on each side simplifies, and we obtain by recalling \eqref{eq:I'a,a}:
$$\sum_{i<j,k<l}K^{ \mathcal{B}_{ij}\mathcal{B}_{kl}}(f\vert m)= -\jap{I_a'(f),\Delta_a f}\geq 2I_a(f).$$
We now follow the proof of Theorem~\ref{thm:main} but use this bound instead of dropping the corresponding term. Pick $\eta=1$ to fix ideas, and write $\alpha=2-\varepsilon$. By inspecting \eqref{eq:derivD4}, we see that for $\chi$ small enough, it holds that
$$\frac{d}{dt}D(f)\leq -2I_a(f\vert m)-(2-\varepsilon/2) I_{\Id}(f\vert m)-2(d-1)\Vert \chi\Vert^2$$
But $D(f)=I_{a+\Id+\chi}(f)-d=I_{a+\Id+\chi}(f\vert m)+(d-1)\Vert \chi\Vert^2$ by Lemma~\ref{lem:relat}, so
as soon as $\chi\leq \varepsilon/(2\alpha)\Id$, it then holds that
    \begin{align*}
        \frac{d}{dt}D(f)&\leq -2I_a(f\vert m)-\alpha I_{\Id }(f\vert m)-(\varepsilon/2)I(f\vert m)-2(d-1)\Vert \chi\Vert^2\\
        &\leq -\alpha I_a(f\vert m)-\alpha I_{\Id }(f\vert m)-\alpha I_\chi(f\vert m)-\alpha(d-1)\Vert \chi\Vert^2=-\alpha D(f).
    \end{align*}
This holds for $\chi$ sufficiently small depending on $\alpha$. Since $\max_i\vert\chi_i\vert\leq e^{-4dt/(d-1)}$, we can indeed find a time $t_1$ depending only on $\alpha$ and $d$ such that for all $t\geq t_1$, $\frac{d}{dt}D(f)\leq -\alpha D(f)$. We conclude by Gronwall's lemma.
\end{proof}
\begin{remark}
    Integrating from $t_1$ to $+\infty$ the above inequality, we recover
    $$\alpha(H(f_{t_1})-H(m))=\alpha\int_{t_1}^\infty D(f_s)ds \leq -\int_{t_1}^\infty \frac{d}{ds}D(f_s)=D(f_{t_1}),$$
    This estimate is known as Cercignani's conjecture, known to hold for the Landau-Maxwell equation, see \cite{DesvillettesVillani2000}. The above provides a proof of the estimate \emph{à la} Stam.
\end{remark}

\begin{remark}
    \label{rem:fullmono}
    We discuss why we expect full monotonicity and how it would require to better understand the ignored terms. Suppose $d=2$, and consider an initial data with $T_1$ very small and $T_2$ close to $2$. This corresponds to $\chi_1 = 1-\varepsilon$ and $\chi_2= -1+\varepsilon$. From \eqref{eq:derivD2}, we see that $\frac{d}{dt}D(f)$ then contains a term like $(12-O(\varepsilon)) I_2(f)$. To absorb it, we can try and apply the same technique that Villani used in \cite{Villani2000} to prove the monotonicity of Fisher information. It consists in writing that for any $\lambda$,
    $$K^{\partial_2 \partial_2}(f) \geq 4\lambda I_2(f) -2\lambda^2$$
    (which is proven by expanding $\int f (\partial^2_{22}\log f + \lambda)^2$ and integrating by parts the cross-term), choosing $\lambda$ large enough to absorb the $I_2$ term, and then dealing with the remaining constants. However the coefficient in front of $K^{\partial_2 \partial_2}(f)$ is $(1+\chi_2)^2= \varepsilon^2$ (instead of $\varepsilon$ for the Fisher information in \cite{Villani2000}). This requires picking $\lambda$ like $\varepsilon^{-2}$ rather than $\varepsilon^{-1}$ which is too large to adapt the sequel of the proof of \cite{Villani2000}. A more promising direction is to try and exploit the other second-order terms, which have a better-behaved coefficient in front of them: in particular, the fully spherical term $K^{ \mathcal{B}_{12}\mathcal{B}_{12}}(f)$ has a constant coefficient in front of it. It vanishes on radial functions, but on the contrary for very anisotropic functions it certainly explodes as $\varepsilon\rightarrow 0$, and we might extract from it a control of $I_2$.
\end{remark}
\sloppy
\printbibliography

\end{document}